\documentclass[12pt,twoside]{article}
\usepackage{amsmath}
\usepackage{hyperref}
\usepackage{amsfonts}
\usepackage{amssymb}

\usepackage{graphics}
\usepackage[pdftex]{graphicx}

\def\bee{\begin{equation}}
\def\eee{\end{equation}}

\input epsf

\begin{document}

\centerline{\Large\bf Criteria  equivalent to  the Riemann Hypothesis}

\bigskip
\bigskip

\centerline{\large {J. Cis{\l}o, M. Wolf}}
\bigskip
\centerline{ Instytut Fizyki Teoretycznej}
\centerline{Uniwersytet Wroc{\l}awski}
\centerline{PL-50-204 Wroc{\l}aw}
\centerline{Pl. Maxa Borna 9}

\bigskip

\begin{abstract}
We give a brief overview of a few criteria equivalent to the Riemann Hypothesis. Next we
concentrate on the Riesz and B{\'a}ez-Duarte criteria. We proof that they are equivalent
and we provide some computer data to support them. It is not compressed to six pages version
of the talk delivered by M.W. during the XXVII Workshop on Geometrical Methods in Physics,
28 June -- 6 July, 2008, Bia{\l}owie{\.z}a, Poland. 
\end{abstract}

\bigskip
\bigskip

\bibliographystyle{plain}
\bigskip

{\bf 1. Introduction}\\

Euler   investigated the series
\bee
\zeta(s)=\sum_{n=1}^\infty \frac{1}{n^s}
\label{zeta_1}
\eee
for  real $s>1$. In particular, he found expression giving values of the zeta function
for even arguments:
\bee
\zeta(2m)= \frac{|B_{2m}|\pi^{2m}}{2(2m)!}
\eee
where $B_n$ are the Bernoulli numbers. Riemann showed \cite{Riemann} that the integral ($s\neq1$):
\bee
\zeta(s)=\frac{\Gamma(-s)}{2\pi i}\int_{+\infty}^{+\infty} \frac{(-x)^s}{e^x-1} \frac{dx}{x}
\eee

\noindent where the integration is performed on the contour

\vskip 3.4cm
\begin{picture}(0,0)(0,0)

\put(160,30){\vector(0,1){60}}
\put(100,60){\vector(1,0){120}}

\put(160,60){\oval(20,20)[l]}
\put(160,65){\oval(10,10)[rt]}
\put(160,55){\oval(10,10)[rb]}

\put(165,55){\vector( 1,0){45}}
\put(210,65){\vector(-1,0){45}}

\end{picture}

\noindent is well  defined on the whole complex plane without $s=1$,
where $\zeta(s)$ has the simple pole, and is equal to (\ref{zeta_1})
on the right of the line $\Re[s]=1$. The $\zeta(s)$ function has trivial zeros:
$s = -2,-4,-6,\dots$ i.e. $\zeta(-2n)=0$.
Besides that there exists the  infinity of zeros $\rho=\sigma + it$  in the
 {\it critical strip}, $0 \leq  \Re [\rho] =\sigma \leq 1 $. If $\rho$ is zero, then
also  $\overline{\rho}$ and $1-\rho$ are zeros, thus zeros are located symmetrically
around the {\it critical line} $\Re[s]=\frac{1}{2}$, see e.g. \cite{Titchmarsh}.
The Riemann Hypothesis (RH)  states that all non-trivial zeros of $\zeta(s)$ lie on
the critical line  $s=\frac{1}{2}+it$. Presently the requirement that they  are simple is
often added. It is one of the best known open problems in mathematics, see e.g.
\cite{Derbyshire},\cite{Edwards}. In the last years, after the Clay Mathematics Institute
granted 1 million US\$ award for solving the dilemma of the RH, there have appeared on
the arxiv
many preprints claiming to have proved  (e.g.,\cite{Shinya-2007},  \cite{LI-2008}) or disproved
(e.g.  \cite{Aizenberg-2007}, \cite{Pati-2007}) the RH, but so far  all of them
have been withdrawn as errors were found in them.

Already Riemann  calculated numerically a few first nontrivial zeros of $\zeta(s)$
\cite{Edwards}. Below  there is a short list of numerical determinations of nontrivial
zeros of $\zeta(s)$:


J.P. Gram(1903): 15 zeros are on the critical line \cite{Gram}

$$ \vdots $$

A.Turing (1953): 1104 zeros are on the critical line \cite{Turing}

$$ \vdots $$
D.H. Lehmer (1956): 25000 zeros are on the critical line.

$$ \vdots $$

 A few years ago S. Wedeniwski (2005)
was leading the  internet project Zetagrid \cite{Zetagrid} which during four years
determined that  $250\times 10^{12}$ zeros are on the critical line:
 $s=\frac{1}{2}+it$, $~~|t| < 29,538,618,432.236$.
The present record belongs to  K. Gourdon(2004) \cite{Gourdon}: the first  $10^{13}$ zeros
are on the critical line.\\

\bigskip

{\bf 2. Some Criteria equivalent to the RH}\\

There are probably  well over one hundred  statements equivalent to the RH, see eg.
\cite{Titchmarsh},
\cite{aimath}, \cite{Watkins}. Riemann's original aim  was to prove  the guess made by
15-years-old Gauss, namely that the number $\pi(x)$ of primes $<x$  ~is well aproximated by
the logarithmic integral:
\bee
\pi(x)\approx {\rm Li}(x) = \int_2^\infty \frac{du}{\log(u)}.
\eee
In this spirit in 1901, Koch  proved \cite{Koch} that the
Riemann Hypothesis is equivalent to the following error term for the  expression for
the prime counting function: $\pi(x)$ is given by:
\bee
\pi(x) = {\rm Li} (x) + {\mathcal{O}}(\sqrt{x}\ln(x)).
\eee

Another similar criterion is
\bee {\rm RH}  \Leftrightarrow~ \pi(x) = {\rm Li} (x) + {\mathcal{O}}(x^{1/2+\epsilon})~~~
{\rm  ~for~ each~~} \epsilon >0.
\eee

The following criterion is of interest to mathematical physicists.
 In  1955  Arne Buerling \cite{Buerling} proved that the
Riemann Hypothesis is equivalent to the assertion that $\mathcal N_{(0,1)} $
is dense in $L^2(0,1).$ Here $\mathcal N_{(0,1)} $ is the space of functions
\bee
N_{(0,1)}=\left\{ \sum_{k=1}^n c_k \rho\left(\frac{\theta_k}{x}\right), ~~0<\theta_k<1,\sum_{k=1}^n c_k=0,~~n=1,2, 3,...~~\right\}
\eee
where $\rho(u)= u - \lfloor u \rfloor$ is a fractional part of  $u$.
The function $\zeta(s)$ does not have zeros in the half-plane
$\sigma>\frac{1}{q},~1<q<\infty$ {\ iff}  the set $N_{(0,1)}$ is dense in $L^q(0,1)$.

In fact Beurling proved that the following three statements regarding a number
 $q \in (1,\infty)$ are equivalent:

(1) $\zeta(s)$ has no zeros in $\sigma>1/q$

(2) $\mathcal N_{(0,1)} $ is dense in $L^q(0,1)$

3) The characteristic function
$\chi_{(0,1)}$ is in the closure of $\mathcal N_{(0,1)} $ in $L^q(0,1)$\\

The following ideas  show that the validity of the RH is very delicate and subtle.
Let us introduce the function
\bee
\xi(iz)=\frac{1}{2}\left(z^2-\frac{1}{4}\right)\pi^{-\frac{z}{2}-\frac{1}{4}}\Gamma\left(\frac{z}{2}+\frac{1}{4}\right)
\zeta\left(z+\frac{1}{2}\right).
\eee
We can see from  the above formula that the RH  $\Leftrightarrow$ all zeros of $\xi(iz)$
are real. The point is that $\xi(z)$ can be expressed as the following Fourier transform:
\bee
\frac{1}{8}\xi\left(\frac{z}{2}\right)=\int_0^\infty \Phi(t)\cos(zt) dt,
\eee
where
\bee
\Phi(t)=\sum_{n=1}^\infty (2\pi^2 n^4 e^{9t} - 3\pi n^2 e^{5t}) e^{-\pi n^2e^{4t}}.
\eee
And now we follow the rule of Polya \cite{Polya} : if one can not solve a particular problem, maybe
it is possible to solve more general problem. So, we introduce the family of functions
$H(z,\lambda)$ parameterized by $\lambda$ as the following Fourier transform:
\bee
H(z, \lambda) =\int_0^\infty \Phi(t)e^{\lambda t^2} \cos(zt) dt.
\eee
Thus we have $H(z,0)=\frac{1}{8} \xi(\frac{1}{2}z)$
N. G. De Bruijn \cite{de_Bruijn} proved that (1950):

1. $H(z,\lambda)$ has only real zeros for $\lambda\geq \frac{1}{2}$

2. If $H(z,\lambda)$ has only real zeros for some  $\lambda'$, then
$H(z,\lambda)$ has only real zeros for each $\lambda>\lambda'$.

\noindent And here comes  bad news: in 1976 Ch. Newman  \cite{Newman} has proved that
there exists parameter $\lambda_1$ such that
$H(z,\lambda_1)$ has at least one non-real zero. Thus, there exists
such constant $\Lambda$ in the interval $-\infty < \Lambda <\frac{1}{2}$
that  $H(z,\lambda)$ has real zeros $\Leftrightarrow \lambda>\Lambda$.
The Riemann  Hypothesis is equivalent to $\Lambda\leq 0$. This constant $\Lambda$ is now
called the de\_Bruijn--Newman constant.
Newman believes that  $\Lambda \geq 0 $. The computer determination has
provided the numerical values of de\_Bruijn--Newman constant, here is a sample of
results:
$$
{\rm Csordas ~{\it et~ al}~(accuracy~ of~ calculations ~\rm 360 ~digits)} ~~(1988)  -50<\Lambda
$$
$$
\vdots
$$
$$
{\rm te~ Riele~ (accuracy~ of~ calculations ~250 ~~digits) ~~(1991)}~~~~~~~~~-5 <\Lambda
$$
$$
\vdots
$$
$$
{\textstyle\rm Odlyzko~~ (2000)}~~-2.7\cdot 10^{-9}<\Lambda
$$
This last result \cite{Odlyzko2000} was obtained using the following pair of zeros of the Riemann
$\zeta(s)$ function near zero number $10^{20}$ that are unusually close together (thus
they almost violate the simplicity of zeros):
\bee
k=10^{20}+71810732, ~~ \gamma_{k+1}-\gamma_k < 0.000145
\eee
Because the gap  in which $\Lambda$ catching the RH is so squeezed, Odlyzko noted
in \cite{Odlyzko2000}, that ``...the Riemann Hypohesis, if true,  is just barely true''.

Li Criterium (1997) \cite{Li-1997}: Riemann Hypothesis is true  {\it iff}
the sequence:
\begin{displaymath}
\lambda_n=\frac{1}{(n-1)!}\frac{d^n}{ds^n}(s^{n-1}\log \xi(s)) |_{s=1} 
\end{displaymath}
where
\begin{displaymath}
\xi(s)=\frac{1}{2}s(s-1)\Gamma\left(\frac{s}{2}\right)\zeta(s)
\end{displaymath}
fulfills:
\bee
\lambda_n\ge 0~~~~ {\rm for}~~  n=1,2,\dots
\label{Li}
\eee
Explicit expression: $~~\lambda_n =\sum_\rho (1-(1-1/\rho)^n)$. K. Ma{\'s}lanka
\cite{Maslanka-2004}, \cite{Maslanka-2004-b}
gave explicit expression for $\lambda_i$ and performed extensive computer
calculations of these constants confirming (\ref{Li}).

\begin{figure}
\begin{minipage}{15.8cm}
\vspace{-0.3cm}
\begin{center}
\hspace{-2cm} \includegraphics[height=15cm,angle=0, bb=0 0 600 600, scale=1]{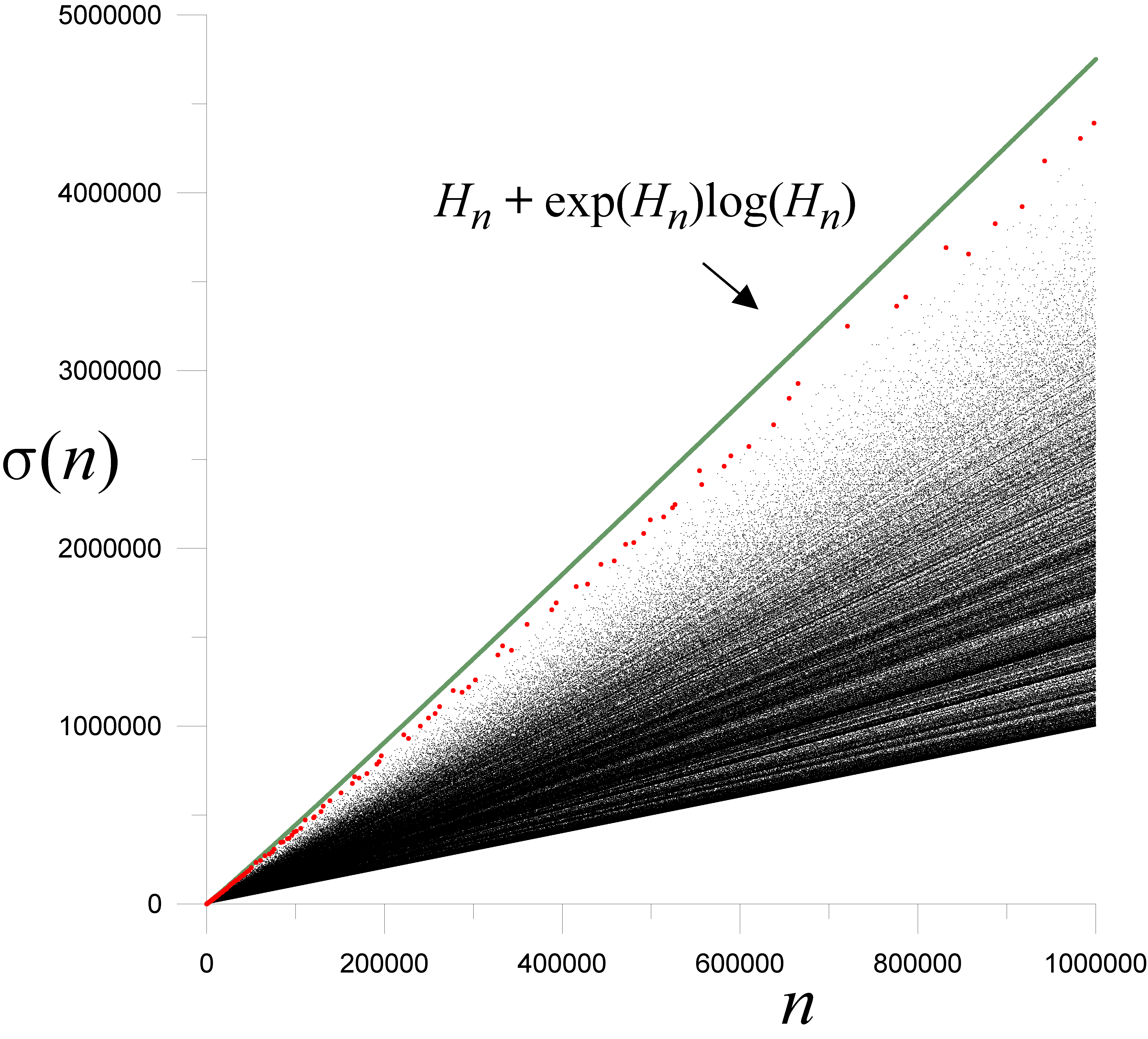} \\
\hspace{-2cm} Fig.1 The plot of $\sigma(n)$ for $1<n<10^6$. In red are plotted values of
  $\sigma(n)$ which  approach the threshold values closer than 5\%. Data for this plot was
obtained  with the free package PARI/GP \cite{Pari}\\
\end{center}
\end{minipage}
\end{figure}

A sensation was stirred  in 2000, when the elementary criterion for the RH was invented 
by Lagarias \cite{Lagarias}:
the Riemann Hypothesis is equivalent to the  inequalities:
\bee
\sigma(n) \equiv \sum_{d|n} d \leq H_n + {\rm{e}}^{H_n}\log(H_n)
\label{Lagarias}
\eee
for each  $n=1,2,\ldots$,  where $H_n$ is the $n$-th harmonic number
$ H_n=\sum_{j=1}^n \frac{1}{j}.$ Here $\sigma(n)$ is the sum of all divisors of $n$.
The illustration of (\ref{Lagarias}) is shown on Fig.1.  In the paper \cite{Briggs} the maxima
of $\sigma(n)$ were studied.\\

\bigskip
\newpage

{\bf 3. Criteria of Riesz and B{\'a}ez-Duarte}\\

In 1916 Riesz \cite{Riesz}  introduced  the function:
\bee
R(x) = \sum_{k=0}^\infty \frac{(-1)^{k}x^{k+1}}{k!\zeta(2k+2)}
\eee
and next he  proved that:
\bee
RH \Leftrightarrow   R(x) = \mathcal{O}\left(
x^{1/4+\epsilon}\right).
\label{Riesz}
\eee
It involves uncountably many values of $x$ and L. B{\'a}ez-Duarte \cite{Baez-Duarte}
 found a discrete version of (\ref{Riesz}). In turn, his derivation was based on the
representation for $\zeta(s)$ found by  K. Ma{\'s}lanka  \cite{Maslanka}:
\bee
\zeta(s)=\frac{1}{(s-1)\Gamma(1-s/2)}\times \sum_{k=0}^\infty \frac{\Gamma(k-\frac{s}{2})}{k!}\sum_{j=0}^k(-1)^j \binom{k}{j}(2j+1)\zeta(2j+2).
\label{M_repr}
\eee

Using (\ref{M_repr}) L. B{\'a}ez-Duarte  proved that the Riemann Hypothesis is equivalent
to the fact that the sequence
\bee
c_k=\sum_{j=0}^k {(-1)^j \binom{k}{j}\frac{1}{\zeta(2j+2)}}
\label{c-k}
\eee
decreases to zero like
\bee
c_k={\mathcal{O}}(k^{-\frac{3}{4}+\epsilon})~~~~~{\rm for~~ each}~~ \epsilon.
\eee
In 2005 one of us \cite{Wolf-2006}   started the computer
calculation of $c_k$. The plot of $c_k$ is presented in Fig.2. The envelops are
given by the equations
\bee
y(k)=\pm A k^{-\frac{3}{4}},~~~~ A=0.777506\ldots \times 10^{-5}.
\eee

\begin{figure}
\begin{minipage}{15.8cm}
\vspace{-0.3cm}
\begin{center}
\hspace*{-13cm}\includegraphics[width=15truecm,angle=270, bb=0 0 600 600, scale=0.08]{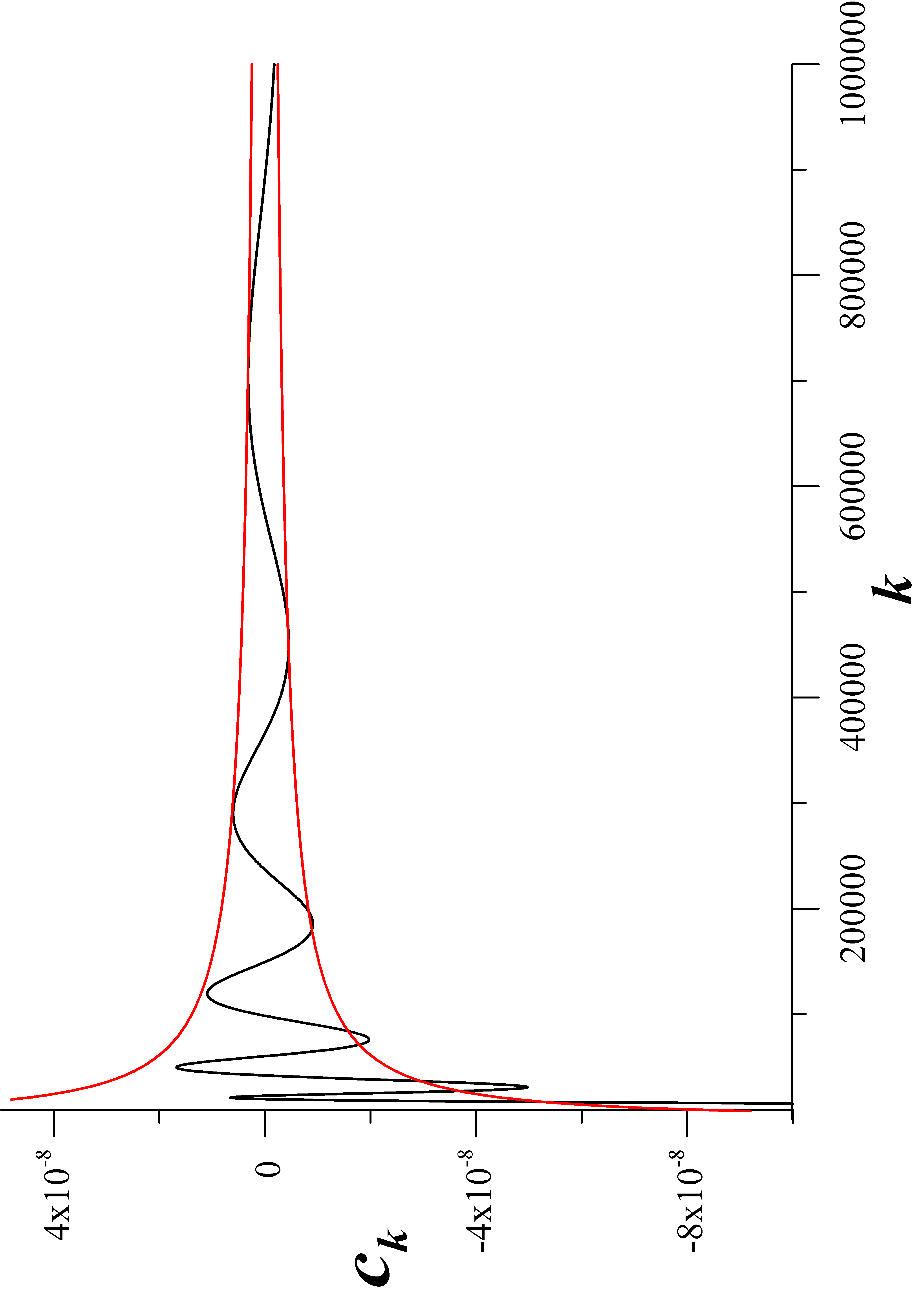} \\
\vspace{8cm}
\hspace{-2cm}\vspace{-0.7cm} Fig.2 The plot of $c_k$ given by (\ref{c-k}) \\
\end{center}
\end{minipage}
\bigskip
\bigskip
\end{figure}

In \cite{Cislo-Wolf-2006} and \cite{Cislo-Wolf-2008}  we have shown that $c_k$ and $R(x)$ are ``entangled'' by means of
the relation:
\bee
\sum_{k=0}^\infty   \frac{c_k x^k}{k!} = \frac {e^x}{x} R(x)
\label{Cislo2}
\eee
Next we have proved
that\\
{\bf Theorem: } For any real number $ \delta > -3/2 $ we have
\begin{equation} \label{Th}
R(x) = O(x^{\delta+1}) \Leftrightarrow  c_k = O(k^\delta).
\end{equation} \\
The proof of this theorem is based on the fact that $R(k)/k\approx c_k$ for integer $k$. In more
detail, in \cite{Cislo-Wolf-2008} we have shown that:\\
\begin{equation}
\left| \frac{R(k)}{k} -c_k \right| \leq \frac{3\sqrt{\pi}}{16} k^{-3/2} + \mathcal{O}(k^{-2}).
\label{bound}
\end{equation}
We have also obtained  the value of the sum:
\bee
\sum_{k=0}^{\infty}(-1)^k c_k  = \sum_{k=1}^{\infty}\frac{1}{2^k} \frac{1}{\zeta(2k)}= 0.7825279853253842\ldots.
\eee
Finally, we mention the approximate relation which follows from (\ref{bound})
entangling in a strange manner values of $\zeta(2j+2)$:
\bee
\sum_{j=0}^\infty \frac{(-1)^{j}k^{j}}{j!\zeta(2j+2)}\approx
\sum_{j=0}^k (-1)^j \binom{k}{j}\frac{1}{\zeta(2j+2)}.
\eee
We have checked numerically that the difference between these two sums very
quickly tends to zero  with increasing $k$.\\

{\bf Acknowledgment:} We thank
Dr Anna Cis{\l}o for the reading and polishing of the manuscript.\\


\end{document}